\documentclass[1pt,notitlepage,twoside,a4paper]{amsart}

\usepackage{amsmath,amssymb,enumerate}

\usepackage{epsfig,fancyhdr,color}

\usepackage{amssymb}
\usepackage{amsmath,amsthm}    
\usepackage{latexsym}
\usepackage{amscd} 
\usepackage{psfrag}
\usepackage{graphicx} 
\usepackage[latin1]{inputenc}  
\usepackage[all]{xy} 
\usepackage[mathcal]{eucal}

\newtheorem{proposition}{\rm\bf Proposition}[section]

\theoremstyle{definition}

\theoremstyle{remark}
\newtheorem{remark}[proposition]{\rm\bf Remark}
\newtheorem{remarks}[proposition]{\rm\bf Remarks}

\def\interieur#1{\mathord{\mathop{\kern 0pt #1}\limits^\circ}}

\definecolor{NoteColor}{rgb}{1,0,0}

\title[On five papers 
 by Herbert Gr\"otzsch]{On five papers 
 by Herbert Gr\"otzsch}
\author{Vincent Alberge and Athanase Papadopoulos}
\address{V. Alberge, Fordham University, Mathematics Department, 441 East Fordham Road, Bronx, NY 10458, U.S.A; A. Papadopoulos, Institut de Recherche Math{\'e}matique Avanc\'ee,
Universit{\'e} de Strasbourg and CNRS,
7 rue Ren\'e Descartes,
 67084 Strasbourg Cedex, France}

 \date{\today}

 \makeindex

\begin{document}

\maketitle
 \begin{abstract} 
Herbert Gr\"otzsch is the main founder of the theory of quasiconformal mappings.  We review five of his papers, written between 1928 and 1932,  that show the progress of his work from conformal to quasiconformal geometry.  This will give an idea of his motivation for introducing quasiconformal mappings, of the problems he addressed and on the results he obtained concerning these mappings.

The final version of this paper will appear in Vol. VII of the \emph{Handbook of Teichm\"uller theory}.
\end{abstract} 
 
   \bigskip
   


\section{Introduction}
    
Herbert Gr\"otzsch is the main founder of the theory of quasiconformal mappings. In a series of papers written between 1928 and 1932, he introduced these mappings as a natural generalization of conformal mappings and he developed their main properties. He saw that several problems in conformal geometry naturally lead to problems on quasiconformal mappings. 
In this chapter, we review five of his papers that show the progress of his work from conformal to quasiconformal geometry. The five papers are the following:

\begin{enumerate}
\item \label{1}  ``\"Uber einige Extremalprobleme der konformen Abbildung" (On some extremal problems of the conformal mapping), published in 1928
 \cite{Groetzsch1928};

\item \label{2}   ``\"Uber  einige Extremalprobleme der konformen Abbildung.  II" (On some extremal problems of the conformal mapping. II), published in 1928
  \cite{Groetzsch1928b};

\item  \label{3}    ``\"Uber  die Verzerrung bei schlichten nichtkonformen Abbildungen und \"uber eine damit zusammenhängende Erweiterung des Picardschen Satzes."
(On the distortion of univalent non-conformal  mappings and a related extension of the Picard theorem) \cite{gr3}; published in 1928.

\item \label{4}    ``\"Uber die Verzerrung bei nichtkonformen schlichten Abbildungen mehrfach zusammenhängender schlichter
Bereiche" (On the distortion of non-conformal schlicht
mappings of multiply-connected schlicht regions), published in 1930 \cite{gr4};

\item \label{5}  ``\"Uber  m\"oglichst konforme Abbildungen von schlichten Bereichen"
  (On closest-to-conformal mappings of schlicht domains), published in 1932 \cite{gr5}.

\end{enumerate}

Translations of these five papers into English are included in the present volume.

We shall present the main results of these papers, providing introductory remarks, establishing the connections between them, and explaining the background. We shall also indicate some relations with works of Lavrentieff and Teichm\"uller.

  The first two papers  are concerned with conformal representations of multiply-connected domains\index{conformal representation!multiply-connected domains}.   Gr\"otzsch\index{Gr\"otzsch, Herbert (1902--1993)} uses there, for the first time,  the so-called \emph{length-area method}, a powerful method which is a direct application of the Schwarz inequality.  In the third paper, he introduces the notion of \emph{non-conformal} (``nichtkonformen") mapping, which he also calls \emph{mapping of bounded infinitesimal distortion} (``Abbildung von beschr\"ankter infinitesimales Verzerrung"). This notion is very close to that of quasiconformal mapping, as it is understood today. The results he proves are all valid for quasiconformal mappings (with the same proofs) and for this reason we shall refer to Gr\"otzsch's non-conformal mappings as quasiconformal mappings. His aim in this paper is to show that several properties known for conformal mappings remain true in this more general setting. He proves a distortion theorem\index{Theorem!distortion} for the new class of mappings---an analogue of a distortion theorem of Koebe---and a generalization of the big Picard\index{Picard theorem (big)}\index{Theorem!Picard (big)}\index{Picard theorem (big)! for quasiconformal mappings}\index{Theorem!Picard (big)! for quasiconformal mappings}  theorem, one of the very classical results on meromorphic functions. The fourth paper is concerned with quasiconformal representations of multiply-connected domains.\index{quasiconformal representation!multiply-connected domains} 
   The results obtained in this paper are generalizations of results known for conformal representations. The fifth paper is concerned with maps that Gr\"otzsch calls  \emph{closest-to-conformal} (``m\"oglichst konforme"). These are the mappings that we call today "extremal quasiconformal mappings".\index{extremal quasiconformal mapping}\index{quasiconformal mapping!extremal}

 \section{Conformal representations}

 The problem of mapping conformally domains of the complex plane $\mathbb{C}$ onto ``canonical" domains was inaugurated by Riemann in his doctoral dissertation \cite{Riemann-Grundlagen} (1851), in which he proved that any simply-connected open subset of the plane which is not the whole plane can be conformally mapped in a one-to-one manner onto the unit disc. He also proved that such a mapping is unique up to  post-composition by a conformal automorphism of the disc (a M\"obius transformation). This is the well-known \emph{Riemann Mapping Theorem}.\index{Riemann Mapping Theorem}\index{Theorem!Riemann Mapping} Riemann proved it using the so-called Dirichlet principle, which he formulated and used extensively in his works on Riemann surfaces. Generalizing the Riemann Mapping Theorem to non-simply-connected planar domains is a broad subject that occupied generations of mathematicians after Riemann. We recall the following two well-known facts:
\begin{enumerate}
\item Two connected open subsets of the complex plane that have the same connectivity are not necessarily conformally equivalent. For instance, two circular annuli (that is, annuli in the complex plane bounded by two concentric circles) are conformally equivalent if and only if they have the same \emph{module},\index{module of an annulus} that is, if and only if the ratio of the two radii is the same for both annuli.

\item There is no natural ``canonical" class of domains for multiply-connected domains. Circular annuli may be considered as some kind of ``standard" domains for doubly-connected regions, but there are other possibilities. \end{enumerate}
We now discuss some of these standard multiply-connected domains.

The unit disc slit along an interval of the form $[0,r]$ ($r<1$) is a useful object in conformal geometry, and it was employed by Gr\"otzsch\index{Gr\"otzsch, Herbert (1902--1993)}  in his work on conformal and quasiconformal geometry. It is known in the classical literature on quasiconformal mappings under the name \emph{Gr\"otzsch domain}\index{Gr\"otzsch extremal domain}\index{extremal domain!Gr\"otzsch}. Such a domain appears in particular as a solution of an extremal quasiconformal problem described in Theorem 9 below, which is Theorem 3 of Gr\"otzsch's  paper \cite{gr4}. The same extremality property of this domain is proved in the book by Lehto and Virtanen, \emph{Quasiconformal 
mappings in the plane}, \cite[p. 52]{LV73}. In the same book, the Gr\"otzch domain is shown to be the solution of the following extremal problem for conformal mappings: if $\mu(r)$ denotes the module of this domain, then,  the module of any doubly-connected domain separating the points $0$ and $r$ from the circle $\vert z\vert =1$ is at most equal to $\mu(r)$; see \cite[p. 54]{LV73}. 
 
 In Ahlfors' book on conformal invariants, the expression \emph{Gr\"otzsch annulus} is used for the complement, in the complex plane, of the closed unit disc and a segment of the form $[R,\infty]$ for some positive $R$ \cite[p. 72]{AhC}. It is possible to write explicitly a conformal homeomorphism between the last two domains.

 Teichm\"uller, in his paper  \emph{Untersuchungen \"uber konforme und quasikonforme Abbildung} (Investigations on conformal and quasiconformal mappings), translated in the present volume,  calls a \emph{Gr\"otzsch extremal region} the complement of the unit disc in the Riemann sphere $\mathbb{C}\cup\{\infty\}$  cut along a segment of the real axis joining a point $P>0$ to the point $\infty$ 
 \cite[\S 2.1]{T13}. He gives estimates for this domain and he uses it in his investigations.

Another ``standard model" for doubly-connected domains is the Riemann sphere $\mathbb{C}\cup\{\infty\}$ slit along  two intervals of the form $[-r_1,0]$ and $[r_2,\infty)$ where $r_1$ and $r_2$ are positive numbers. This domain is known under the name \emph{Teichm\"uller extremal domain}\index{Teichm\"uller extremal domain}\index{extremal domain!Teichm\"uller}  \cite[p. 52]{LV73}.  Teichm\"uller uses this domain in his study of extremal properties of topological annuli, in \S 2 of his paper \cite{T13}; see in particular \S 2.6. This domain satisfies the following extremal property: its module is an upper bound for the module of any doubly-connected domain that separates the pair $0$ and $-r_1<0$ from the pair $r_2>0$ and $\infty$. Teichm\"uller proved this characterization in  \cite[\S 2.4]{T13}, using Koebe's one-quarter theorem and Koebe's distortion theorem.\index{Theorem!Koebe!distortion}\index{Theorem!Koebe!one-quarter}\index{Koebe distortion theorem}\index{Koebe one-quarter theorem} His proof is reproduced in Ahlfors' book \cite[p. 72]{AhC}; see also \cite[p. 55]{LV73}. 
 The module of the Teichm\"uller extremal domain is equal to the quantity $2\mu\sqrt{\frac{r_1}{r_1+r_2}}$ (where $\mu$ denotes as above the module of the Gr\"otzsch extremal domain).
 
  Teichm\"uller, in his paper 
\cite{T13},  also works with other ``standard" domains, e.g.  the circular annulus $1<\vert z\vert <P_2$ cut along the segment joining $z=P_1$ to $z=P_2$, where $P_1$ is a point on the real axis satisfying $1<P_1<P_2$ \cite[\S 2.4]{T13}. He proves several extremal properties of these domains. He uses them in a geometric proof which simplifies and generalizes a distortion result that Ahlfors obtained in his thesis \cite{ahlforsthesis}; see also the discussion in Ahlfors' book \cite[p. 76]{AhC}.

Among the other ``standard" doubly-connected domains, we mention the \emph{Mori domain},\index{Mori extremal domain}\index{extremal domain!Mori} called so in the monograph \cite{LV73}.


\begin{figure}[!h] 
\centering
\psfrag{-r1}{\small $-r_1$}
 \psfrag{0}{\small $0$}
  \psfrag{r2}{\small $r_2$}
   \psfrag{infty}{\small $\infty$}
    \psfrag{1}{\small $1$}
     \psfrag{r}{\small $r$}
\includegraphics[width=1\linewidth]{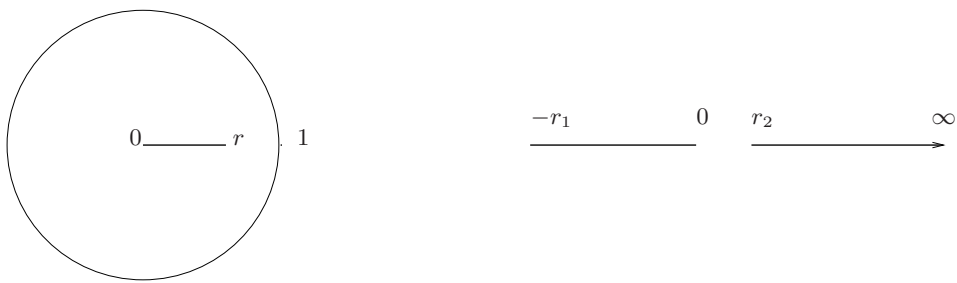} 
\caption{\smaller A Gr\"otzsch extremal domain (left) and a Teichm\"uller extremal domain (right), called so in the book of Lehto and Virtanen \cite{LV73}. In the figure to the right, the doubly-connected domain is the complement in the Riemann sphere of the union of the segment and the ray indicated.}   
\label{GT-domains}  
\end{figure}

Before Gr\"otzsch, Koebe\index{Koebe, Paul (1882--1945)} studied extensively conformal mappings of finitely-connected domains of the plane onto \emph{circle domains},\index{circle domain} that is, multiply-connected domains in the plane whose boundary components are all circles (which may be reduced to points), cf. \cite{Kobe1, Kobe2}. He proved that every finitely-connected domain in the plane is conformally equivalent to a circle domain. This generalizes the Riemann Mapping Theorem, which says that a simply-connected domain which is not the entire plane is conformally equivalent to the unit disc. We mention by the way that the question of whether there is a result analogous to Koebe's theorem which is valid for an open subset of the plane having an infinite number of boundary components is open, and is known as the \emph{Kreisnormierungsproblem},\index{Koebe!Kreisnormierungsproblem} or the \emph{Koebe uniformization conjecture}.\index{Koebe!uniformization conjecture} Koebe formulated this conjecture in his 1908 paper \cite{Kobe1908} p. 358. We refer the reader to the recent paper \cite{Bowers} by Bowers in which the author surveys many developments of this conjecture, including the works on circle packings done by Koebe, Thurston and others.

\begin{figure}[!h] 
\centering
\includegraphics[width=.8\linewidth]{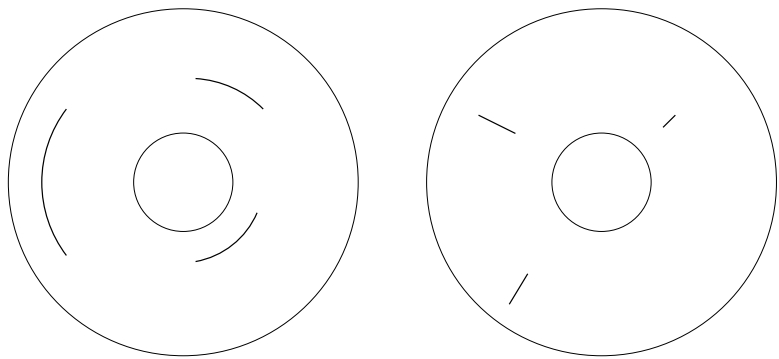} 
\caption{\smaller An annulus with circular slits (left) and an annulus with radial slits (right)}   
\label{Fig-circ-slit}  
\end{figure}

Gr\"otzsch, who was a student of Koebe, was naturally led to study problems connected with conformal representations. In the paper \cite{Groetzsch1928}, he works with two classes of domains which were already considered by Koebe:
\begin{enumerate}
 \item  \emph{Annuli with circular slits}.\index{annulus with circular slits} These are circular annuli from which a certain number ($\geq 0$) of circular arcs (which may be reduced to points) centered at the center of the annulus, have been removed, see Figure \ref{Fig-circ-slit} (left). 
 \item \emph{Annuli with radial slits}.\index{annulus with radial slits}  These are circular annuli from which a certain number ($\geq 0$) of radial arcs (which may be reduced to points) have been removed. A radial arc is a segment in the annulus which, when extended, passes through the center of the annulus, see Figure \ref{Fig-circ-slit} (right). 
 \end{enumerate}
 
 Koebe, in the paper \cite{Groetzsch1928}, showed that any finitely-connected domain of the complex plane which is not the whole plane admits a conformal representation\index{Koebe!conformal representation}\index{conformal representation!Koebe} (which is essentially unique) onto a domain of one of the two kinds above.

 \section{Two lemmas from the paper \cite{Groetzsch1928}}
 
Gr\"otzsch's paper \emph{ \"Uber einige Extremalprobleme der konformen Abbildung} \cite{Groetzsch1928} starts with two lemmas that give a subadditivity result for the modules of (finitely or infinitely many) disjoint subdomains of a circular annulus in the complex plane.\index{subadditivity for modules of annuli}
In the first lemma, the subdomains are disjoint topological quadrilaterals, each having two opposite sides on the boundary circles of the circular annulus. In the second lemma, the subdomains are disjoint topological annuli that are homotopy-equivalent to the circular annulus.

In each case, the lemma says that the sum of the modules of the subdomains is bounded above by the module of the ambient circular annulus.  In particular, if we divide the interior of a circular annulus into two annuli by a simple closed curve homotopic to the boundary components, then the sum of the modules of the two resulting annuli is not greater than the modules of the ambient annulus. 

Gr\"otzsch also obtains a characterization of the equality case: equality occurs, in the first case, if and only if the rectangles are obtained by radial slicing of the circular domain, and in the second case, if and only if the decomposition into subdomains is obtained by a slicing of the circular domain by circles concentric to the boundary. Furthermore, in the two equality cases, the subdomains must fill out the ambient circular annulus.

The proofs of these module inequalities\index{module inequality (annului)} are based on the so-called \emph{length-area} method.\index{length-area method}
  Regarding this method, Ahlfors\index{Ahlfors, Lars Valerian (1907--1996)}  writes, in  his collected papers edition \cite{Ahlfors-Collected}, commenting on his first two published papers on the asymptotic values of entire functions of finite order, which in fact are those of his doctoral dissertation: ``[\ldots] The early history of this method is obscure, but I knew it from and was inspired by its application in the well-known textbook of Hurwitz-Courant\footnote{Courant and Hurwitz, in their book \emph{Funktionentheorie} (Springer, Berlin, 1922), used the length-area method if their proofs of results of Carath\'eodory on the boundary values of the Riemann mapping, which the latter has published in then three papers \cite{Cara1, Cara2, Cara3}.} to the boundary correspondence in conformal mapping." Talking about Nevanlinna and himself, he adds: ``None of us was aware that only months earlier H. Gr\"otzsch had published two important papers on extremal problems in conformal mapping in which the same method is used in a more sophisticated manner [\ldots]  The method that Gr\"otzsch and I used is a precursor of the method of extremal length." In fact, Ahlfors applied this method  in his thesis \cite{ahlforsthesis} in order to prove what became later known as the \emph{Ahlfors distortion theorem}.
 Teichm\"uller, in his papers \cite{T13,T20} and \cite{T23}, uses extensively  the length-area method, and he calls it the \emph{Gr\"otzsch--Ahlfors} method.  

The module inequalities proved in \cite{Groetzsch1928} are considered today as classical results. They were also obtained (with different methods) and widely generalized by Teichm\"uller in his paper \emph{Untersuchungen \"uber konforme und quasikonforme Abbildung} (Investigations on conformal and quasiconformal mappings) (1938), translated in the present volume (\S 2.4 and 2.5 of \cite{T13}).

\section{Four theorems  from the paper \cite{Groetzsch1928}}

In the paper \cite{Groetzsch1928}, Gr\"otzsch uses the two lemmas of the previous section in his proof of four theorems which we now state. Theorems 3 and 4 answer questions posed by S. Szeg\"o, see \cite{S1}.

For $0<r<1$,  we denote by $\mathcal{K}_r$ the annulus $r<\vert z\vert<1$ and by $\overline{\mathcal{K}_r}$ an annulus with circular slits obtained from $\mathcal{K}_r$  by removing a finite number of circular arcs centered at the origin (Figure \ref{Fig-circ-slit},  left).

 \medskip

\noindent {\bf Theorem 1 (\cite{Groetzsch1928} \S 2).--} \emph{Let $R$ be a real number satisfying $R\geq r$. If a domain $\overline{\mathcal{K}_r}$ is mapped conformally and bijectively onto a domain having $\vert z\vert =r$ and $\vert z\vert =R$ as boundary components with the circle $\vert z\vert =r$ sent onto the circle 
$\vert w\vert =r$ and the circle $\vert z\vert =1$ sent onto the circle  $|w|=R$, then $R\geq 1$. Furthermore, the case $R=1$ holds if and only if the conformal mapping is a rotation.}

\medskip
Theorem 1 implies in particular that if a circular annulus with circular slits and with boundary circles $\vert z\vert=r$ and $\vert z\vert=1$ is mapped conformally and bijectively onto a domain in the complex plane that has the same circles $\vert z\vert=r$ and $\vert z\vert=1$  as boundary components, then the second domain is also an annulus with circular slits and the conformal map is  a rotation.

\medskip

In the next two theorems, a domain $\mathcal{K}_r$ or $\overline{\mathcal{K}_r}$ is  mapped conformally and bijectively by a mapping $w=f(z)$ onto a domain called $\mathcal{B}_r$ (respectively $\overline{\mathcal{B}_r}$)  of the extended complex plane $\mathbb{C}\cup\{\infty\}$ not containing $\infty$ in its interior, such that the circle $|z|=r$ is sent to the circle $|w|=r$ and such that the image of $\mathcal{K}_r$ (respectively $\overline{\mathcal{K}_r}$) is contained in the subset $|w|\geq r$ of $\mathbb{C}\cup\{\infty\}$. 

The image by $w$  of the circle $\vert z\vert =1$ is called the \emph{outer boundary} of $\mathcal{K}_r$ (respectively $\overline{\mathcal{K}_r}$).

Let $d$ be the shortest distance from a point on the outer boundary of $\mathcal{B}_r$ (respectively $\overline{\mathcal{B}_r}$) to the point $w=0$.\index{extremal problem!closet-to-boundary point}

 \medskip
\noindent {\bf Theorem 2 (\cite{Groetzsch1928} \S 3).--} \emph{Consider a conformal mapping of $\mathcal{K}_r$ sending $|z|=r$ to $|w|=r$ and $|z|=1$ to a ray whose extension contains the point $w=0$ and joining a point $q$ to $\infty$. Then,  $d\geq \vert q\vert$.}

\medskip 

\begin{figure}[!h] 
\centering
\includegraphics[width=0.8\linewidth]{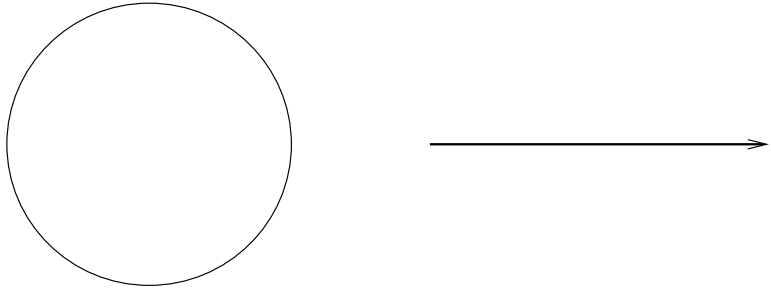} 
\caption{\smaller An annulus with circular boundary with the other boundary being an infinite ray. This extremal domain appears in the second theorem}   
\label{fig-ray}  
\end{figure}

In Theorem 2, the domain represented in Figure \ref{fig-ray} appears as an extremal domain. This theorem is also known under the name  \emph{Gr\"otzsch module theorem}, and it is usually stated under a form which uses the modulus of an annulus instead of the distance.

 \medskip
 In the next theorem, we use the above notation $\mathcal{B}_r$ and $\overline{\mathcal{B}_r}$. 
\medskip

\noindent {\bf Theorem 3 (\cite{Groetzsch1928} \S 3).--} 
\emph{Assume there exist $n$ points on the outer boundary of the image domain $\mathcal{B}_r$ (respectively $\overline{\mathcal{B}_r}$) realizing the distance of this outer boundary to the point $w=0$, assume that these points are the vertices of a regular $n$-gon having $w=0$ as center and that there exists a positive constant $M$ (possibly equal to $+\infty$) such that  $|w(z)|\leq M$ for all $z$. 
Then a map that realizes the shortest distance of the outer boundary of the image domain to the point $w=0$ sends $\mathcal{K}_r$ 
 or $\overline{\mathcal{K}_r}$
 respectively to a domain bounded by $|z|=r$, $|z|=M$ and which has $n$ radial slits whose extensions contain $0$ and which join the points 
$$\alpha\cdot \rho_{r,M,n} \cdot e^{2\pi i\,k\over n}\quad (r<\rho_{r,M,n}<M\,,\ |\alpha|=1\,,k=0,1,\dots,n-1)$$
to points on the circle $|w|=M$. Furthermore, the mappings satisfying the stated properties are the unique ones that realize the required minimum property.}

\begin{figure}[!h]
\centering
\includegraphics[width=0.4\linewidth]{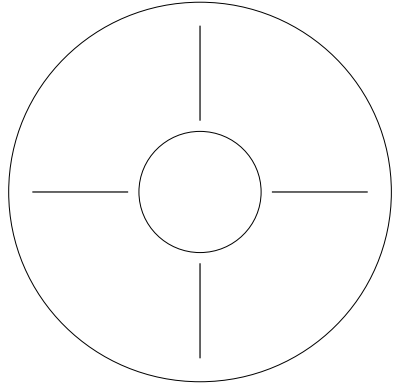} 
\caption{\smaller An extremal domain for Theorem 3, with $k=4$. The inner slits admit a rotational symmetry}   
\label{fig-symmetric}  
\end{figure}

\medskip
The extremal region that appears in Theorem 3 is represented in Figure \ref{fig-symmetric}.\index{extremal problem!closet-to-boundary point}
\medskip

In the next theorem, $\mathcal{K}_0$ and $\overline{\mathcal{K}_0}$ denote respectively the unit disc in $\mathbb{C}$ and the unit disc slit along finitely many circular arcs centered at the origin. When we map the domain $\mathcal{K}_0$ or $\overline{\mathcal{K}_0}$ by a conformal map $f$, we shall say, in analogy with the previous cases, that the image by $w$ of  the circle $\vert z\vert =1$ is the \emph{outer boundary} of the image  $f(\mathcal{K}_0)$ or $f(\overline{\mathcal{K}_0})$.

\medskip

\begin{figure}[!h]
\centering
\includegraphics[width=0.5\linewidth]{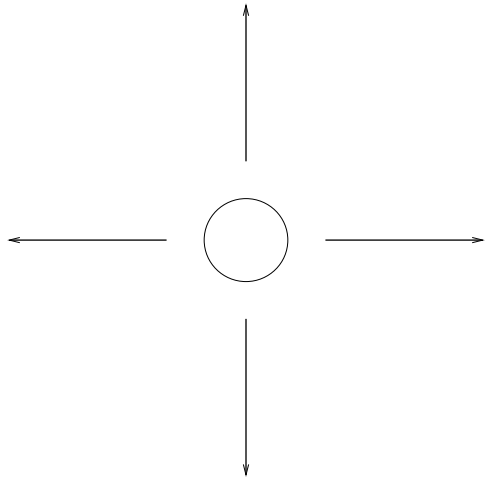} 
\caption{\smaller An extremal domain for Theorem 4. The inner slits admit a rotational symmetry}   
\label{fig-symmetry2}  
\end{figure}

\noindent {\bf Theorem 4 (\cite{Groetzsch1928} \S 4).--}
\emph{Suppose that $\mathcal{K}_0$  (respectively $\overline{\mathcal{K}_0}$) is mapped conformally and bijectively by $w=f(z)$ onto a domain $\mathcal{B}_0$ (respectively $\overline{\mathcal{B}_0}$) not containing $\infty$ in its interior satisfying $f(0)=0$ and $|f'(0)|=1$. Suppose furthermore that there exist $n$ points on the outer boundary of $\mathcal{B}_0$ (respectively $\overline{\mathcal{B}_0}$) that are the vertices of a regular $n$-gon centered at $0$ and whose distance to the point $w=0$ is the shortest distance from a point on the outer boundary to that point. Then 
$$d\geq \root n\of {1\over 4}\,,$$
and this extremal value is attained by the mapping $f(z)={z\over \root n\of {(1+z^n)^2}}$. This mapping sends $\mathcal{B}_0$ to the $w$--plane slit along the $n$ rays emanating from the points
$\root n\of{1\over 4}\cdot e^{2\pi i\,k\over n}$ ($k=0,1,2,\dots,n-1$) to $\infty$ (Figure \ref{fig-symmetry2}.)}
 
\medskip

\medskip

\begin{remarks}
1.--- An extremal result similar to the one that appears in Theorem 2 was obtained by Lavrentieff in his paper \cite{L13} (1934). Lavrentieff attributes the result to Gr\"otzsch, and he uses completely different methods (analysis instead of geometry).

2.--- Extremal problems similar to those studied by Gr\"otzsch that we mentioned in this section are analyzed by Teichm\"uller in his paper \cite{T13} (1938), in particular \S 2.2 to 2.3. Teichm\"uller attributes the results he obtains to Gr\"otzsch.

3.--- 

\end{remarks}
  \section{Two theorems from the paper \cite{Groetzsch1928b}}
 
The paper  (\emph{On Extremal Problems for conformal mappings, II}) \cite{Groetzsch1928b} 
 by Gr\"otzsch is a sequel to the paper \cite{Groetzsch1928} and  is concerned with the same subject.\index{extremal problem!conformal mappings}
The proofs are also based on the lemmas proved in the first paper.\index{conformal representation}

The notation is as follows (we use Gr\"otzsch's notation so that the reader can easily compare with the paper):

  The unit  disc in the complex plane is denoted by $\mathfrak{R}_0$. 
  
  For $r>0$,  $\mathfrak{R}_r$ is the  subset of the complex plane defined by $r\leq\vert z\vert <1$.

A conformal mapping $w=f(z)$  defined on $\mathfrak{R}_0$ is said to be \emph{normalized} if it satisfies $f(0)=0$ and $\vert f'(0)\vert =1$, and if $\infty$ is not in the image.

A conformal mapping $w=f(z)$ defined on $\mathfrak{R}_r$ is said to be \emph{normalized} if the circle $\vert z\vert =r$ is sent to itself  and if $\infty$ is not in the image.

Gr\"otzsch establishes distortion theorems for the modules $\vert f(z)\vert$ and $\vert f'(z)\vert$ of normalized injective conformal mappings  defined on $\mathfrak{R}_r$ ($r\geq 0$). The formulation of the result uses a conformal mapping $E_r(z)$ defined as follows:

For $r>0$, $E_r(z)$ is the conformal map that maps $\mathfrak{R}_r$ onto a domain of the $w$-plane bounded by the circle $\vert w\vert =r$ and a suitable slit contained in the positive real axis going to infinity and whose extension in the finite direction contains the point $+r$ on the real line.

$E_r(z)$ is the conformal map that maps $H_r$ onto the domain of the .... 

Cette application est unique!

We set 
\[m_r(\rho)= \min_{\vert z\vert =\rho} \vert E_r(z)\vert\]
and
\[M_r(\rho)= \sup_{\vert z\vert =\rho} \vert E_r(z)\vert.\]

\medskip

\noindent {\bf Theorem 5.--} \emph{We have, for all $\vert z\vert =\rho$: 
\[m_r(\rho)\leq \vert f(z)\vert\leq M_r(\rho)\]
with the two inequalities being equalities holding if and only if $f(z)=\alpha E_r(\beta z)$
with $\vert\alpha\vert=\vert\beta\vert=1$.}
\medskip

For $\mathfrak{R}_0$, a similar statement hods, and in this case an explicit formula for the mapping $E_0$ and the extremal values $m_0(\rho)$ and  $M_0(\rho)$ are given.
In fact, the resulting map is the so-called \emph{Koebe map}, and the values of $m_0(\rho)$ and  $M_0(\rho)$ are those that are given by Koebe's distorsion map.

In the next theorem, Gr\"otzsch gives estimates on the derivative $\vert f'(z)\vert$. Here, the notation is as follows:

\[m'_r(\rho)= \min_{\vert z\vert =\rho} \vert E'_r(z)\vert\]
and
\[M'_r(\rho)= \sup_{\vert z\vert =\rho} \vert E'_r(z)\vert.\]

The result is then the following:

\medskip

\noindent {\bf Theorem 6.--} We have
\[m'_r(\rho)\leq \vert f'(z)\vert\leq M'_r(\rho)\]
for $\vert z\vert =\rho$, with equality holding only for the functions $\alpha E_r(\beta z)$ as above.

\medskip
In the case of $\mathfrak{R}_0$, explicit formulae for $m'_r(\rho)$ and $M'_r(\rho)$ are given.

These results generalize Koebe's distortion theorems.

\section{On the content of Gr\"otzsch's paper \cite{gr3}}
The paper \cite{gr3} is the first which contains Gr\"otzsch's notion of quasiconformal mappings. This is a mapping with bounded infinitesimal distortion $Q$.  He introduces such a class of mappings in the first section of the paper and he denotes them by $\mathfrak{A}_{Q}$. As a matter of fact, the term \emph{nichtkonformen} (which we translate by non-conformal) is only used in the title. It was Ahlfors, in \cite{Ahlfors-Zur} who used for the first time the term \emph{quasikonform} (``quasiconoformal'');\index{quasiconformal!mapping} cf.\  Ahlfors' comments at the beginning of Volume 1 of his \emph{Collected Papers} \cite{Ahlfors-Collected}.

For $Q\geq 1$, Gr\"otzsch defines in \cite{gr3} a mapping in the class $\mathfrak{A}_Q$ to be a bijective mapping of a domain $G$ of the $z$-plane onto a surface which covers the Riemann sphere or a domain in it (Gr\"otzsch  uses the term ``spread as a Riemann surface over a $w$-plane,'' which is close to the terminology used by Riemann) which satisfies the following three properties:
\begin{enumerate}
\item $f$ is continuous on $G$;
\item up to a set with at most countably many interior  points (called exceptional points) of $G$, the mapping $f$ is a local diffeomorphism which sends an infinitesimal circle to an infinitesimal ellipse whose ratio of minor to major axis is between $\frac{1}{Q}$ and $Q$;
\item  In the neighborhood of each exceptional point, the map $f$ is a finite branched covering.
\end{enumerate}
Such a  mapping from a domain of the Riemann sphere onto another domain of the Riemann sphere is not necessarily  bijective, but it can be lifted to a bijective map from a domain of the Riemann sphere onto a Riemann surface which is a ramified covering of the sphere. Furthermore,  the mapping is not assumed to be sense-preserving. 
 Lehto and Virtanen call these ``nicht-conformal maps'' of \cite{gr3} ``regular quasiconformal" (see \cite[p.17]{LV73}). 

In the second section of the paper, Gr\"otzsch establishes the property known today as the \emph{geometric definition}\index{quasiconformal mapping!geometric definition} of quasiconformality, saying that the image of a quadrilateral of module $M$ by a bijective mapping in $\mathfrak{A}_Q$ is a quadrilateral whose module $M^{\prime}$ satisfies the inequality
\begin{equation}\label{eqgrine}
\frac{1}{Q} M \leq M^{\prime} \leq Q M.
\end{equation} 
 
  The proof that Gr\"otzsch gives of this property is based on the length-area method, which he already used in the papers \cite{Groetzsch1928} and \cite{Groetzsch1928b}.

 After setting the inequality (\ref{eqgrine}), Gr\"otzsch writes that in the case of rectangles,   equality in (\ref{eqgrine}) holds ``only for certain immediately determinable extremal affine mappings.'' This result is now called the \emph{solution of the Gr\"otzsch Problem}; cf. Ahlfors' book \cite[p. 8]{ahlforsqc}). The solution was given by Gr\"otzsch in the article \cite{gr5} which we review below. 
 
  Right after the relation (\ref{eqgrine}), Gr\"otzsch deduces the corresponding inequalities for annuli. More precisely, using the logarithm function and then the exponential function, he obtains the fact that if there exists a bijective mapping in the class $\mathfrak{A}_Q$ from the annulus with inner radius $r$ and outer radius $1$ onto the annulus with inner radius $\widetilde{r}$ and outer radius $1$, then
\begin{equation}\label{eqgrine2}
r^{Q}  \leq \widetilde{r} \leq r^{\frac{1}{Q}}.
\end{equation} 
As in the case of rectangles, Gr\"otzsch says that equalities hold only for certain simple extremal maps which consist of rotations and reflections. By using logarithms, we can deduce from the double inequality (\ref{eqgrine2}) the geometrical definition  of quasiconformal mappings for doubly-connected domains.

In the same section of the paper, Gr\"otzsch notes that a result of Carathéodory  is still valid for bijective $\mathfrak{A}_Q$ mappings, namely, any bijective $\mathfrak{A}_Q$ mapping from a simply-connected domain bounded by a Jordan curve onto a domain of the plane can be extended continuously to a homeomorphism. 

  To conclude this section, Gr\"otzsch states two distortion inequalities\index{distortion inequality!annuli} that hold for an arbitrary mapping $f$ in the class $\mathfrak{A}_Q$  defined on the annulus $0\leq r<\left| z \right| < 1$, sending the circle of radius $r$ onto the circle of radius $r$ and such that  for any $r\leq \left| z \right| < 1$, we have $\left| f(z) \right| \geq r$. The first inequality concerns the behaviour of the distance between the circle of radius $r$ and a point in the image by $f$ of any circle of radius $\rho$ with $r < \rho < 1$. The second inequality concerns the behaviour of the rate of change of $f$.

In the last section, Gr\"otzsch gives an extension of the (big) Picard theorem for mappings in the class $\mathfrak{A}_Q$.\index{Picard theorem (big)}\index{Theorem!Picard (big)}\index{Picard theorem (big)! for quasiconformal mappings}\index{Theorem!Picard (big)! for quasiconformal mappings} More precisely, he proves that for a given $r$, any $\mathfrak{A}_Q$ mapping defined on the punctured disc $0<\left| z \right|<r$ for which $0$ is an essential singularity and which omits at most two points in the plane (or three points in the Riemann sphere) is constant. (An essential singularity is a point at which the function does not admit any finite or infinite limit.)

The idea of Gr\"otzsch's proof is quite elementary. Indeed, he takes  an $\mathfrak{A}_Q$ mapping $f$ defined on the punctured disc $0<\left| z \right| <r$, extending continuously to the set $\left| z \right| = r$ with an essential singularity at $0$ and which omits three values in the Riemann sphere, that we can assume to be $0$, $1$ and $\infty$. By definition of an $\mathfrak{A}_Q$ map, we can lift $f$ to a bijective map $\widetilde{f}$ from $0<\left| z \right|<r$ onto a Riemann surface $\mathcal{S}_f$ equipped with a holomorphic map onto a subset of the Riemann sphere such that $f = \pi \circ \widetilde{f}$.  
Since $\mathcal{S}_f$ is doubly-connected, it is biholomorphic by a mapping $\varphi$ to an annulus $a < \left| z \right| < b$. Since $0$ is an essential singularity, we have $a>0$, because otherwise the map $\pi \circ \varphi^{-1}$ would be a non-constant  holomorphic map from $0 < \left| z \right| < b$ to $\mathbb{C}\setminus\left\lbrace 0, 1\right\rbrace$ with an essential singularity at $0$, which is impossible by the big Picard theorem (for meromorphic functions). Thus, we have an $\mathfrak{A}_Q$ mapping $\varphi \circ \widetilde{f}$ from $0<\left| z \right| <r$ onto $a < \left| z \right| < b$, which is impossible because of Relation (\ref{eqgrine2}) and the supperadditivity for the module of an annulus.

In conclusion, Gr\"otzsch in this short paper gives a certain number of results which show that conformal and quasiconformal mappings share several properties.

\section{Some comments on Gr\"otzsch's paper \cite{gr3}}

 1.-- Picard's big theorem\index{Picard theorem (big)}\index{Theorem!Picard (big)}\index{Picard theorem (big)! for quasiconformal mappings}\index{Theorem!Picard (big)! for quasiconformal mappings} says that if a holomorphic function has an essential singularity, then in any punctured neighborhood of this essential singularity, the function, takes infinitely often all possible values in the Riemann sphere, with at most two exceptions. 
 Even though the statement that refers to Picard in Gr\"otzsch's paper is not an exact analogue of that theorem, 
  what is proved gives easily such a result. Indeed, if the quasiconformal mapping omits at most two points in the plane,  we can lift it to the associated Riemann surface which is  necessarily a punctured disc (since otherwise we would have a quasiconformal mapping onto an annulus), and thus, we get a holomorphic function with one essential singularity. Therefore (by Picard's theorem) the quasiconformal mapping takes all values infinitely often with the exception of at most two.

\medskip

2.-- The idea of a class of functions which generalize conformal mappings is also the subject of Lavrentieff's paper \cite{lavrentieff} in which this author gives another proof of the big Picard theorem for such maps. 

\medskip

3.-- The use of supperadditivity in \S 2 of Gr\"otzsch's paper, along with the double inequality (\ref{eqgrine2}), was used by Teichm\"uller in \cite{T13}   in order to prove that quasiconformal mappings preserve the type of a simply connected Riemann surfaces.

\medskip

 4.-- Gr\"otzsch continued to investigate similarities between conformal and quasiconformal mappings in his paper \cite{gr4} which we review next, by solving geometrical extremal problems for quasiconformal mappings, after he studied such problems in the setting of conformal mappings in the papers \cite{Groetzsch1928, Groetzsch1928b}.

\section{The paper \cite{gr4}}

Gr\"otzsch's paper \emph{\"Uber die Verzerrung bei nichtkonformen schlichten Abbildungen mehrfach
zusammenh\"angender schlichter Bereiche} \cite{gr4}  
  published in 1930  is a sequel to the two papers \cite{Groetzsch1928} and  \cite{Groetzsch1928b}. Here, instead of working with conformal representations, Gr\"otzsch considers representations of multi-connected domains by quasiconformal mappings,\index{quasiconformal mapping} which he calls  ``nichtkonformen".  The terminology is close to te one he uses in his paper \cite\cite{gr3}. Here, Gr\"otzsch means by this word a mapping defined on a region $B$ of the $z$-plane which can be approximated uniformly, in the neighborhood of each point except  for finitely many points in $B$---called exceptional points---by an affine mapping, and such that in the neighborhood of each non-exceptional point the ratio $a/b$ of the great axis to the small axis of the infinitesimal ellipse which is the image of an infinitesimal circle is uniformly bounded by two constants: 
 \[\frac{1}{Q}\leq  a/b \leq Q
\]
where $Q$ is independent of the choice of the non-exceptional point. 

As we did for the mappings introduced by Gr\"otzsch in his paper \cite{gr3}, we shall use here the terminology $\mathfrak{A}_Q$ mapping for such a mapping.

Gr\"otzsch's results in this paper are based again on the length-area method, which he used extensively in his previous papers in the setting of conformal mappings, together with the  two double inequalities (\ref{eqgrine}) and (\ref{eqgrine2}) which he proved in his paper \cite{gr3}. The method, as in the case of conformal mappings, consists in taking appropriate sequences of surface strips.

Gr\"otzsch's goal is to determine, under some normalization conditions, the mappings in the class $\mathfrak{A}_Q$  that realize extremal values to certain distortion quantities.
He obtains six theorems which we state below. Before stating them, we introduce some notation (which is Gr\"otzsch's notation).

$B_{n+1}$ is an $(n+1)$-connected open region of the complex plane. It is said to be \emph{normalized} if the circle $\vert z\vert =1$ is one of its boundary components and if it is contained in the interior of that circle. 

All the regions $B_{n+1}$ that we consider are normalized. We say that $\vert z\vert =1$ is the \emph{outer boundary} of $B_{n+1}$, and we let $R_1, \ldots, R_{n}$ be its remaining boundary components  (they may be reduced to points).

 Gr\"otzsch gives the solutions of eight extremal theorems. The question, in each case, is to find a mapping in the class $\mathfrak{A}_Q$ of normalized mappings for which some geometrically defined quantity is extremal.
 
  For each of these problems, any $\mathfrak{A}_Q$ mapping which is a solution turns out to have constant dilatation $Q$.
  
In the first 4 theorems, a \emph{normalized} mapping  $\mathfrak{A}_Q$ of $B_{n+1}$ is a bijective mapping from $B_{n+1}$ onto a normalized region $\tilde{B}_{n+1}$ in the complex plane sending the outer boundary $\vert z\vert =1$ of  $B_{n+1}$ to the outer boundary $\vert w\vert =1$ of  $\tilde{B}_{n+1}$.
We let $\tilde{R}_1, \ldots, \tilde{R}_{n}$ be the images of the boundary components $R_1, \ldots, R_{n}$.

In the first theorem, the quantity that is maximized is the area of the domain enclosed by $\tilde{R}_k$.
  
  \medskip
  \noindent {\bf Theorem 7 (Theorem 1 of \cite{gr4}).--} In the set $\mathfrak{A}_Q$ of normalized mappings of $B_{n+1}$, the area of the domain enclosed by $\tilde{R}_k$ attains its maximum if and only if 
  \begin{enumerate}
  \item $\tilde{R}_k$ is a circle centered at $w=0$;
  \item the other  $\tilde{R}_p$ (if they exist) are circular arc slits centered at $w=0$;
  \item at each point of the image, the minor axis of all the distortion ellipse (which Gr\"otzsch calls the Tissot indicatrix)\index{Tissot indicatrix} is directed radially toward the origin $w=0$;
  \item the distortion is everywhere constant and equal to $Q$.
  \end{enumerate}
  Furthermore, any two such extremal mappings differ from each other by at most a rotation around $w=0$ and a reflection in a straight line through $w=0$.
  
  \medskip
  
  The proof of this theorem, like the one of the theorems that follow, is first given in the case of conformal mappings ($Q=1$), and then for $\mathfrak{A}_Q$ mappings, using the two distortion results for quadrilaterals and annuli that we recalled.

    The remaining seven theorems are of the same sort. Instead of  restating them, we indicate, for each one, the quantities that are minimized. In each case, Gr\"otzsch gives a description of the extremal image domain. 
        
  \medskip
  \noindent {\bf Theorem 8  (Theorem 2 of \cite{gr4}).--} In the set $\mathfrak{A}_Q$ of normalized mappings of $B_{n+1}$, to find the mappings for which the maximum value of the diameter of $\tilde{R}_k$  is attained.
  
In this case, the extremal image domain has a rectilinear slit passing through the origin $w=0$ and bisected by that point.

  \medskip
  \noindent {\bf Theorem 9  (Theorem 3 of \cite{gr4}).--} In the set $\mathfrak{A}_Q$ of normalized mappings $w=f(z)$ of $B_{n+1}$ and for a given  pair $(z_1,z_2)$ in $B_{n+1}$, to find the mappings for which the maximum of the quantity $\vert f(z_1)-f(z_2)\vert$ is attained.

  In this case, in the extremal image domain, the images $f(z_1)$ and $f(z_2)$ are the endpoints of a rectilinear segment bisected by the origin $w=0$.
\medskip

  \noindent {\bf Theorem 10  (Theorem 4 of \cite{gr4}).--} Assume that  $B_{n+1}$ has $z=0$ as interior point. In the set $\mathfrak{A}_Q$ of normalized mappings  $w=f(z)$ of $B_{n+1}$ that satisfy $f(0)=0$, and for a given $z_1$ in $B_{n+1}$, to find the mappings for which the maximum of the quantity $\vert f(z_1)\vert$ is attained.

\medskip

In the remaining 4 theorems, a different normalization of the mappings in  $\mathfrak{A}_Q$ is used. Here, an $(n+1)$-connected region $B_{n+1}$ is said to be normalized if it is 
bounded (``from inside") by a circle $ | z | = r $ and lies entirely outside $ | z | = r $, and which is bounded (``from outside") by a linear boundary $L$, either separating $ B_{n + 1}$ from $ \infty $ (which means the curve is closed) or heading towards to $ \infty $,  and  by  $ n-1 $ other boundaries $ R_k $ lying between $ | z | = r $ and $L$.
 
Here, a mapping $\mathfrak{A}_Q$ of $ B_{n + 1}$ is called normalized if the image region $ \tilde{B}_{n + 1}$ of the $w$-plane is also normalized, if  the circle $ | z | = r $ is sent to $ | w | = r $, and if  $L$ is sent to the linear boundary $\tilde{L}$ joining $\tilde{B}_{n + 1}$ to $ \infty $. (In particular, $\infty$ is not an image point).

Let $ \tilde{R}_k $ be  the boundary of  $B_{n + 1}$  corresponding to $ R_k $. The next three theorems are stated in the form of problems of which Gr\"otzsch gives the solution.

  \noindent {\bf Theorem 11  (Theorem 5 of \cite{gr4}).--} Among the normalized mappings of $B_{n+1}$, to find those $\mathfrak{A}_Q$ mappings for which the distance of the boundary $\tilde{L}$ from the point $w=0$ attains its minimum value.
\medskip

  \noindent {\bf Theorem 12  (Theorem 6 of \cite{gr4}).--} Among the normalized mappings of $B_{n+1}$, to find those $\mathfrak{A}_Q$  mappings for which the distance of the boundary $\tilde{L}$ from the point $w=0$ attains its maximum value.
\medskip

  \noindent {\bf Theorem 13  (Theorem 7 of \cite{gr4}).--} Among the normalized mappings $w=f(z)$ of $B_{n+1}$, and for a given point $z_1$ in $B_{n+1}$, to find those $\mathfrak{A}_Q$ mappings for which the minimum of the quantity $\vert f(z_1)\vert$ is attained.
\medskip

  \noindent {\bf Theorem 14  (Theorem 8 of \cite{gr4}).--} Among the normalized mappings $w=f(z)$ of $B_{n+1}$, and for a given point $z_1$ in $B_{n+1}$, to find those $\mathfrak{A}_Q$ mappings for which the maximum of the quantity $\vert f(z_1)\vert$ is attained.
\medskip

\section{The paper \cite{gr5}}

The paper \emph{\"Uber  m\"oglichst konforme Abbildungen von schlichten Bereichen} \cite{gr5}
 by Gr\"otzsch
was published in 1932 and is concerned with maps he calls  \emph{m\"oglichst konforme} (closest-to-conformal).\index{closest-to-conformal mapping} These are maps that we call today ``extremal quasiconformal mappings". Gr\"otzsch starts by recalling that (by the Riemann Mapping Theorem)\index{closest-to-conformal mapping!simply connected domains}\index{extremal mapping!simply connected domains} one can map conformally and in a one-to-one way an arbitrary simply-connected open subset $B_1$ of the plane which is not the whole plane onto another such subset $\tilde{B}_1$, and that this map is completely determined if one asks that the image of three arbitrarily chosen distinct points on the boundary of $B_1$ are sent to three arbitrarily chosen distinct points on the boundary of $\tilde{B}_1$. He then remarks that if one takes, on the boundary of each domain, four cyclically ordered points instead of three, then, generally speaking, one cannot find a conformal mapping between the two domains tat send distinguished points to distinguished points. In this case, one looks for mappings whose deviation from conformality is the smallest. The question is then to introduce an appropriate measure for this deviation from  conformality and to investigate the existence of such closest-to-conformal mappings.

The question is also relevant for domains of the plane which are not simply-connected. In the case of multiply-connected domains, the same question appears naturally, with or without distinguished points on the boundary. In the rest of the paper, Gr\"otzsch gives several examples of closest-to-conformal mappings\index{closest-to-conformal mapping!simply connected domains}\index{extremal mapping!simply connected domains} between domains of the plane, some of them simply-connected and others multiply-connected.

Right at the beginning of the investigation, Gr\"otzsch says that one has to allow the maps to have singularities at isolated points in the interior, possibly countably many and converging to points on the boundary. At the non-singular points (which he calls ``differential-geometric" and which we call ``regular"), the maps are assumed to be differentiable, and the defect in conformality is measured by the ratios of the great axis to the small axis of infinitesimal ellipses which are images of infinitesimal circles, as in the previous papers, \cite{gr3} and \cite{gr4}. Likewise, the quantity $Q\geq 1$ is defined to be the upper bound of such ratios, taken over all the regular points. The measure of deviation of the mapping is taken to be the quantity $Q-1$. Gr\"otzsch notes that this quantity is invariant by (pre- and post-) composition of the mapping with conformal mappings. He then remarks that in a neighborhood of a regular point, infinitesimal squares are sent to infinitesimal parallelograms, and he gives a result on the distortion of sufficiently small parallelograms in the form of a double inequality bounding the quotient of the area over the square of the length of one of the sides in terms of the upper bound $Q$ of the mapping's dilatation. More precisely, the  result says the following (we use Gr\"otsch's notation; formula (1) in his paper):
 
\medskip

\noindent {\bf Proposition.--}
\emph{For every $\epsilon >0$ and at every regular point, there exists a neighborhood of this point such that the image of an infinitesimal square is an infinitesimal parallelogram satisfying the following: if $\tilde s$ is the length of any one of its sides, and $\tilde I$ its area, then
 \begin{equation}\label{star}
(1-\varepsilon)\cdot\dfrac 1Q\cdot \tilde{s}^2<\tilde I<(1+\varepsilon)\cdot Q\cdot \tilde{s}^2,
\end{equation}
where $\varepsilon >0$ can be arbitrarily small.}
\medskip

After this proposition, Gr\"otzsch gives examples of closest-to-conformal mappings.\index{closest-to-conformal mapping!rectangles}\index{extremal mapping!rectangles} He calls the problem of finding such mappings ``the mapping problem,"\index{mapping problem (Gr\"otzsch)} a term that suggests the Riemann Mapping Theorem.

\section{Examples of closest-to-conformal mappings (extremal quasiconformal mappings)}\index{extremal quasiconformal mapping from Gr\"otzsch's paper \cite{gr5}}

\subsection*{Example 1: Rectangles} The mapping problem is considered between two Euclidean rectangles in the plane, where the distinguished points are the vertices.  The rectangles may be assumed to have the same height, and their lengths are denoted by $p$ and $\tilde{p}$ respectively, with $p\leq \tilde{p}$. Gr\"otzsch proves that we always have the inequality $Q\geq \frac{\tilde{p}}{p}$, with equality holding if and only if the map is affine, in which case the value of $Q$ is equal, at every point, to the ratio of the two axes of the image infinitesimal ellipse.

 In his proof, Gr\"otzsch uses the \emph{strip method},\index{strip method}\index{length-area method} a method that is a form of the length-area method. It is also based on  inequalities that relate lengths of curves in a certain family and the area of the domain that the family occupies.
 For this method, he refers to his paper \emph{\"Uber die Verzerrung bei nichtkonformen schlichten Abbildungen mehrfach
zusammenh\"angender schlichter Bereiche}  \cite{gr4}  (Gr\"otzsch calls it so in Example 3, where he uses it again).

\subsection*{Example 2: Doubly connected domains} Gr\"otzsch remarks that any annulus can be mapped conformally onto a Euclidean annulus bounded by two circles. By using logarithms, the problem of finding closest-to-conformal mappings\index{closest-to-conformal mapping!annuli}\index{extremal mapping!annuli} between Euclidean annuli is reduced to that of finding closest-to-conformal mappings between Euclidean rectangles, and hence to Example 1. The unique closest-to-conformal mapping is therefore, in this case too, the affine map. Gr\"otzsch notes that in the case where the boundary component of one annulus (and not the other one) is reduced to a point, then we cannot find a closest-to-conformal mapping between the two (the quantity denoted by $Q$, that is, the upper bound of the major to the minor axes of infinitesimal ellipses is always infinite).

\subsection*{Example 3: Simply-connected domains with two distinguished interior points} Gr\"otzsch considers the problem of finding closest-to-conformal mappings between two simply-connected domains\index{closest-to-conformal mapping!simply connected domains}\index{extremal mapping!simply connected domains} with two distinguished points in the interior. Using the Riemann Mapping Theorem, one may assume that each domain is the open unit disc in the complex plane and that in each such disc, the distinguished points lie on the real axis and are symmetric with respect to the origin. By cutting each of the two discs along the segment joining the two distinguished points, we obtain annuli, and the mapping problem is reduced to that of Example 2. Gr\"otzsch notes that the closest-to-conformal mapping obtained is $C^1$ at the two distinguished points only in the case where it is conformal.

\subsection*{Example 4: simply-connected domains with two distinguished points on the boundary and one distinguished point in the interior.} Gr\"otzsch considers the problem of finding closest-to-conformal mappings between two simply-connected domains\index{closest-to-conformal mapping!simply connected domains}\index{extremal mapping!simply connected domains} having two distinguished points on the boundary and one distinguished point in the interior. Like in the Example 3, by composing with a conformal mapping, each of the two domains can be assumed to be the unit disc in $\mathbb{C}$, with the interior distinguished point being the center of the disc and the distinguished points on the boundary being the points $e^{\pm i\alpha}$ and $e^{\pm i\beta}$ in the first and the second disc respectively, where $\alpha$ and $\beta$ are real numbers satisfying $0<\alpha<\pi$ and $0<\beta<\pi$. We may also  assume that the mapping sought for (supposing it is orientation-preserving) sends the arc of circle whose extremities are the two points $e^{\pm i\alpha}$ and containing the point $-1$ to the arc  of circle containing the two points $e^{\pm i\beta}$ and  containing the point $-1$. One then cuts the two domains along the segment joining $0$ to $-1$, obtaining in each case a simply-connected domain with four distinguished points on the boundary. We are then reduced to finding the closest-to-conformal map between two simply-connected domains with four boundary points on the boundary, and we proceed as in Example 1. The closest-to-conformal map obtained between the two initial discs is $C^1$ except at the origin, and it is $C^1$ there only if it is conformal.

\begin{remark}
A few years after Gr\"otzsch's papers appeared, Teichm\"uller considered the existence and uniqueness problems of closest-to-conformal mappings mappings (which were already called extremal quasiconformal mappings) for arbitrary surfaces: compact or not, orientable or not, with an arbitrary number of distinguished interior or boundary points.  See the paper \cite{T20} \emph{Extremale quasikonforme Abbildungen und quadratische Differentiale} (1939) and the commentary in \cite{T20C}.
In 1944, he published a paper, \emph{Ein Verschiebungssatz der quasikonformen Abbildung} (A displacement theorem for quasiconformal
mapping) \cite{T-displacement},  in which he considers a new type of existence problem for extremal mappings, in which each point on the boundary is a distinguished point. In other words, he studies the question finding and describing the extremal quasiconformal mapping from the unit disc with the origin as a distinguished point to the unit disc with some distinguished point in the interior which is the identity on the boundary (that is, each point on the unit circle is also considered as a distinguished point).
We refer the reader to the paper \cite{VA} where the developments of this problem are also discussed. 
\end{remark}

In the last section of the paper \cite{gr5}, called \emph{Additional remark}, Gr\"otzsch establishes an extension of Schwarz's lemma to the class of differentiable quasiconformal mappings, that is, to the mappings  that belong to the  class denoted by $\mathfrak{A}_{Q}$ in the paper \cite{gr3} discussed above. The extension is valid for branched coverings of the unit disc.

The proof of this extension involves Riemann surfaces which cover the sphere and it uses the classical Schwarz lemma for conformal mappings.\index{Schwarz lemma!quasiconformal mapping}

  This result is in the trend of previous results of Gr\"otzsch whose aim is the extension to the class of quasiconformal mappings of results that are known to hold for conformal mappings; see e.g. the extension of Picard's theorem in the paper \cite{gr3} we reviewed above.

In a biographical note on Gr\"otzsch published in the present volume \cite{Kuehnau-Groetzsch}, Reiner K\"uhnau writes: ``It is strange that today many people associate the name Gr\"otzsch only to the ridiculous `Gr\"otzsch-ring'." We hope that after publishing here the English translations of these five papers of Gr\"otzsch, his work on the theory of quasiconformal mappings will be better known.

 \end{document}